\documentclass{amsart}
\usepackage{amssymb}
\usepackage{amsmath}
\usepackage{amsfonts}

\newtheorem{lemma}{Lemma}[section]
\theoremstyle{definition}
\newtheorem{definition}[lemma]{Definition}
\newtheorem{example}[lemma]{Example}

\theoremstyle{remark}
\newtheorem{remark}[lemma]{Remark}
\numberwithin{equation}{section}
\theoremstyle{plain}

\newtheorem{corollary}[lemma]{Corollary}

\newtheorem{proposition}[lemma]{Proposition}

\input{tcilatex}

\begin{document}
\title[Background notes on Trees]{Some Notes on Trees and Paths}
\author{B.M. Hambly}
\address[B.M. Hambly and Terry J. Lyons]{Mathematical Institute\\
Oxford University\\
24-29 St. Giles\\
Oxford OX1 3LB\\
England}
\email[B.M. Hambly]{hambly@maths.ox.ac.uk}
\author{Terry J. Lyons}
\email[Terry J. Lyons]{tlyons@maths.ox.ac.uk}
\thanks{The authors gratefully acknowledge EPSRC support: GR/R29628/01,
GR/S18526/01}

\begin{abstract}
These notes cover background material on trees which are used in the paper 
\cite{HL}.
\end{abstract}

\maketitle

\section{Trees and paths - background information}

In the paper \cite{HL} it is shown that trees have an important role as the
negligible sets of control theory, quite analogous to the null sets of
Lebesgue integration. The trees considered are \emph{analytic} objects in
flavour, and not the finite combinatorial objects of undergraduate courses. 
%They are easy to define and work with; however, we
%have found that many analysts find the characterisation we use unintuitive.
In this note we collect together a few related ways of looking at them, and
prove a basic characterisation generalising the concept of height function.

We first recall that\newline
(1). Graphs $\left( E,V\right) $ that are acyclic and connected are
generally called \emph{trees}. If such a tree is non-empty and has a
distinguished vertex $\mathbf{v}$ it is called a \emph{rooted tree}.\newline
(2). A rooted tree induces and is characterised by a \emph{partial order on }%
$V$\emph{\ with least element }$\mathbf{v}$. The partial order is defined as
follows 
\begin{equation*}
a\preceq b\quad\text{if the circuit free path from the root } \mathbf{v}%
\rightarrow b\text{ goes through }a.
\end{equation*}
This order has the property that for each fixed $b$ \emph{the set }$\left\{
a\preceq b\right\} $\emph{\ is totally ordered} by $\preceq$.

Conversely any partial order on a finite set $V$ with a least element $v$
and the property that for each $b$ the set $\left\{ a\preceq b\right\} $ is
totally ordered defines a unique rooted tree on $V$. One of the simplest
ways to construct a tree is to consider a (finite) collection $\Omega $ of
paths in a graph sharing a fixed initial or starting vertex, and with the
partial order that $\omega \preceq \omega ^{\prime }$ iff $\omega $ is an
initial segment of $\omega ^{\prime }.$\newline
(3). Alternatively, let $\left( E,V\right) $ be a graph extended into a
continuum by assigning a length to each edge. Let $d\left( a,b\right) $ be
the infimum of the lengths of paths\footnote{%
the sum of the lengths of the edges} between the two vertices $a,\;b$ in the
graph. Then $g$ is a geodesic metric on $V$. Trees are exactly \emph{the
graphs that give rise to 0-hyperbolic metrics} in the sense of Gromov (see
for example \cite{Kapovich95}).\newline
(4). There are many ways to enumerate the edges and nodes of a finite rooted
tree. One way is to think of a family tree recording the descendants of a
single individual (the root). Start with the root. At the root, if all
children have been visited stop, at any other node, if all the children have
been visited, move up to the parent. If there are children who have not been
visited, then visit the oldest unvisited child. At each time $n$ the
enumeration either moves up an edge or down an edge - each edge is visited
exactly twice. Let $h\left( n\right) $ denote the distance from the top of
the family tree after $n$ steps in this enumeration with the convention that 
$h\left( 0\right) =0$, then $h$ is similar to the path of a random walk,
moving up or down one unit at each step, except that it is positive and
returns to zero exactly as many times as there are edges coming from the
root. Hence $h\left( 2\left\vert E\right\vert \right) =0$.

\emph{The function }$h$\emph{\ completely describes the rooted tree.} The
function $h$ directly yields the nearest neighbour metric on the tree. If $h$
is a function such that $h\left( 0\right) =0$, it moves up or down one unit
at each step, is positive and $h\left( 2\left| E\right| \right) =0$, then $d$
defined by 
\begin{equation*}
d\left( m,n\right) =h\left( m\right) +h\left( n\right) -2\inf _{u\in\left[
m,n\right] }h\left( u\right),
\end{equation*}
is a pseudo-metric on $\left[ 0,2\left| V\right| \right] $. If we identify
points in $\left[ 0,2\left| V\right| \right] $ that are zero distance apart
and join by edges the equivalence classes of points that are distance one
apart, then one recovers an equivalent rooted tree.

Put less pedantically, let the enumeration be $a$ at step $n$ and $b$ at
step $m$ and define 
\begin{equation*}
d\left( a,b\right) =h\left( m\right) +h\left( n\right) -2\inf_{u\in \left[
m,n\right] }h\left( u\right) ,
\end{equation*}
then it is simple to check that $d$ is well defined and is a metric on
vertices making the set of vertices a tree.

Thus excursions of simple (random) walks are a convenient (and well studied)
way to describe abstract graphical trees. This particular choice for \emph{%
coding a tree with a positive function on the interval} can be extended to
describe continuous trees. This approach was used by Le Gall \cite{snake} in
his development of the Brownian snake associated to the measure valued
Dawson-Watanabe process. %It is probably well known in
%other contexts - but we have not found it.
%\end{itemize}

\section{$\mathbf{R}$-trees are coded by continuous functions}

One of the early examples of a continuous tree is the evolution of a
continuous time stochastic process, where, as is customary in probability
theory, one identifies the evolution of two trajectories until the first
time they separate. (This idea dates back at least to Kolmogorov and his
introduction of filtrations). Another popular and equivalent approach to
continuous trees is through $\mathbf{R}$-trees (\cite{rtree} p425 and the
references there).

Interestingly, analysts and probabilists have generally rejected the
abstract tree as too wild an object, and usually add extra structure,
essentially a second topology or Borel structure on the tree that comes from
thinking of the tree as a family of paths in a space which also has some
topology. This approach is critical to the arguments used in   \cite{HL}
where tree-like paths are approximated by with simpler tree-like paths in $1$%
-variation. (They would never converge in the `hyperbolic' metric). In
contrast, group theorists and low dimensional topologists have made a great
deal of progress by studying specific symmetry groups of these trees and do
not seem to find their hugeness too problematic.

Our goal in this subsection of the appendix is to prove the simple
representation: that the general $\mathbf{R}$-tree arises from identifying
the contours of a continuous function on a locally connected and connected
space. The height functions we considered on $\left[ 0,T\right] $ are a
special case.

\begin{definition}
An $\mathbf{R}$-tree is a uniquely arcwise connected metric space, in which
the arc between two points is isometric to an interval.
\end{definition}

Such a space is locally connected, for let $B_{x}$ be the set of points a
distance at most $1/n$ from $x$. If $z\in B_{x}$, then the arc connecting $x$
with $z$ is isometrically embedded, and hence is contained in $B_{x}$. Hence 
$B_{x}$ is the union of connected sets with non-empty common intersection
(they contain $x$) and is connected. The sets $B_{x}$ form a basis for the
topology induced by the metric. Observe that if two arcs meet at two points,
then the uniqueness assertion ensures that they coincide on the interval in
between.

Fix some point $v$ as the `root' and let $x$ and $y$ be two points in the $%
\mathbf{R}$-tree. The arcs from $x$ and $y$ to $v$ have a maximal interval
in common starting at $v$ and terminating at some $v_{1}$, after that time
they never meet again. One arc between them is the join of the arcs from $x$
to $v_{1}$ to $y$ (and hence it is the arc and a geodesic between them).
Hence 
\begin{equation*}
d\left( x,y\right) =d\left( x,v\right) +d\left( y,v\right) -2d\left(
v,v_{1}\right) .
\end{equation*}

\begin{example}
Consider the space $\Omega $ of continuous paths $X_{t}\in E$ where each
path is defined on an interval $\left[ 0,\xi \left( \omega \right) \right) $
and has a left limit at $\left[ 0,\xi \left( \omega \right) \right) $.
Suppose that if $X\in \Omega $ is defined on $\left[ 0,\xi \right) $, then $%
X|_{\left[ 0,s\right) }\in \Omega $ for every $s$ less than $\xi $. Define 
\begin{equation*}
d\left( \omega ,\omega ^{\prime }\right) =\xi \left( \omega \right) +\xi
\left( \omega ^{\prime }\right) -2\sup \left\{ t<\min \left( \xi \left(
\omega \right) ,\xi \left( \omega ^{\prime }\right) \right) |\;\omega \left(
s\right) =\omega ^{\prime }\left( s\right) \;\forall s\leq t\right\} .
\end{equation*}
Then $\left( \Omega ,d\right) $ is an $R$-tree.
\end{example}

We now give a way of constructing $\mathbf{R}$-trees. %\begin{lemma}
%If $I$ is a connected and locally connected topological space, and $%
%h:I\rightarrow\mathbb{R}$ is a positive continuous function that attains its
%lower bound, then together they define a pseudo-metric $d$ on $I$ that is
%continuous as a function on $I\times I$ and so that the $d$-equivalence
%classes form a $\mathbf{R}$-tree $T$ and the quotient map onto $T$ is
%continuous.
%\end{lemma}
The basic idea for this is quite easy, but the core of the argument lies in
the detail so we proceed carefully in stages.

Let $I$ be a connected and locally connected topological space, and $%
h:I\rightarrow\mathbb{R}$ be a positive continuous function that attains its
lower bound at a point $v\in I$.

\begin{definition}
For each $x\in I$ and $\lambda\leq h\left( x\right) $ define $C_{x,\lambda}$
to be the maximal connected subset of $\left\{ y\;|\;h\left( y\right)
\geq\lambda\right\} $ containing $x$.
\end{definition}

\begin{lemma}
\label{lem:cxyexist} The sets $C_{x,\lambda}$ exist, and are closed.
Moreover, if $C_{x,\lambda }\cap C_{x^{\prime},\lambda^{\prime}}\neq\phi$
and $\lambda\leq\lambda ^{\prime}$, then 
\begin{equation*}
C_{x^{\prime},\lambda^{\prime}}\subset C_{x,\lambda}.
\end{equation*}
\end{lemma}

\begin{proof}
An arbitrary union of connected sets with non-empty intersection is
connected, taking the union of all connected subsets of $\left\{
y\;|\;h\left( y\right) \geq\lambda\right\} $ containing $x$ constructs the
unique maximal connected subset. Since $h$ is continuous the closure $%
D_{x,\lambda}$ of $C_{x,\lambda}$ is also a subset of $\left\{ y\;|\;h\left(
y\right) \geq\lambda\right\}$. The closure of a connected set is always
connected hence $D_{x,\lambda}$ is also connected. It follows from the fact
that $C_{x,\lambda}$ is maximal that $C_{x,\lambda}=D_{x,\lambda}$ and so is
a closed set.

If $C_{x,\lambda}\cap C_{x^{\prime},\lambda^{\prime}}\neq\phi$ and $%
\lambda\leq\lambda^{\prime}$, then 
\begin{equation*}
x\in C_{x,\lambda}\cup C_{x^{\prime},\lambda^{\prime}}\subset\left\{
y\;|\;h\left( y\right) \geq\lambda\right\},
\end{equation*}
and since $C_{x,\lambda}\cap C_{x^{\prime},\lambda^{\prime}}\neq\phi$, the
set $C_{x,\lambda}\cup C_{x^{\prime},\lambda^{\prime}}$ is connected. Hence
maximality ensures $C_{x,\lambda}=C_{x,\lambda}\cup
C_{x^{\prime},\lambda^{\prime}}$ and hence $C_{x^{\prime},\lambda^{\prime}}%
\subset C_{x,\lambda}$.
\end{proof}

\begin{corollary}
Either $C_{x,\lambda}$ equals $C_{x^{\prime},\lambda}$ or it is disjoint
from it.
\end{corollary}

\begin{proof}
If they are not disjoint, then the previous Lemma can be applied twice to
prove that $C_{x^{\prime},\lambda}\subset C_{x,\lambda}$ and $%
C_{x,\lambda}\subset C_{x^{\prime},\lambda}.$
\end{proof}

\begin{corollary}
If $C_{x,\lambda}=C_{x^{\prime},\lambda}$, then $C_{x,\lambda^{\prime%
\prime}}=C_{x^{\prime},\lambda^{\prime\prime}}$ for all $\lambda^{\prime%
\prime }<\lambda$.
\end{corollary}

\begin{proof}
The set $C_{x,\lambda},C_{x^{\prime},\lambda}$ are nonempty and have
nontrivial intersection. $C_{x,\lambda}\subset C_{x,\lambda^{\prime\prime}}$
and $C_{x^{\prime},\lambda}\subset C_{x^{\prime},\lambda^{\prime\prime}}$
hence $C_{x,\lambda^{\prime\prime}}$ and $C_{x^{\prime},\lambda^{\prime%
\prime }}$ have nontrivial intersection. Hence they are equal.
\end{proof}

\begin{corollary}
$y\in C_{x,\lambda}$ if and only if $C_{y,h\left( y\right) }\subset
C_{x,\lambda}$.
\end{corollary}

\begin{proof}
Suppose that $y\in C_{x,\lambda}$, then $C_{y,h\left( y\right) }$ and $%
C_{x,\lambda}$ are not disjoint. It follows from the definition of $%
C_{x,\lambda}$ and $y\in C_{x,\lambda}$ that $h\left( y\right) \geq\lambda$.
By Lemma~\ref{lem:cxyexist} $C_{y,h\left( y\right) }\subset C_{x,\lambda}$.
Suppose that $C_{y,h\left( y\right) }\subset C_{x,\lambda}$, since $y\in
C_{y,h\left( y\right) }$ it is obvious that $y\in C_{x,\lambda}$.
\end{proof}

%In particular either $C_{x,\lambda}$ equals $C_{x^{\prime},\lambda}$ or it
%is disjoint from it.

\begin{definition}
The set $C_{x}:=C_{x,h\left( x\right) }$ is commonly referred to as the 
\emph{contour} of $h$ through $x$.
\end{definition}

The map $x\rightarrow C_{x}$ induces a partial order on $I$ with $x\preceq y$
if $C_{x}\supseteq C_{y}$. If $h$ attains its lower bound at $x$, then $%
C_{x}=I$ since $\left\{ y\;|\;h\left( y\right) \geq h\left( x\right)
\right\} =I$ and $I$ is connected by hypothesis. Hence the root $v\preceq y$
for all $y\in I$.

\begin{lemma}
Suppose that $\lambda\in\left[ h\left( v\right) ,h\left( x\right) \right]$,
then there is a $y$ in $C_{x,\lambda}$ such that $h\left( y\right) =\lambda$
and, in particular, there is always a contour ($C_{x,\lambda }$) at height $%
\lambda$ through $y$ that contains $x$.
\end{lemma}

\begin{proof}
By the definition of $C_{x,\lambda}$ it is the maximal connected subset of $%
h\geq\lambda$ containing $x$; assume the hypothesis that there is no $y$ in $%
C_{x,\lambda}$ with $h\left( y\right) =\lambda$ so that it is contained in $%
h>\lambda,$ hence C$_{x,\lambda}$ is a maximal connected subset of $%
h>\lambda $. Now $h>\lambda$ is open and locally connected, hence its
maximal connected subsets of $h>\lambda$ are open and $C_{x,\lambda}$ is
open. However it is also closed, which contradicts the connectedness of the $%
I$. Thus we have established the existence of the point $y$.
\end{proof}

The contour is obviously unique, although $y$ is in general not. If we
consider the equivalence classes $x\symbol{126}y$ if $x\preceq y$ and $%
y\preceq x$, then we see that the equivalence classes $\left[ y\right] _{%
\symbol{126}}$ of $y\preceq x$ are totally ordered and in one to one
correspondence with points in the interval $\left[ h\left( v\right) ,h\left(
x\right) \right]$.

\begin{lemma}
\label{lem:loccon} If $z\in C_{y,\lambda}$ and $h\left( z\right) >\lambda$,
then $z$ is in the interior of $C_{y,\lambda}$. If $C_{x^{\prime},\lambda^{%
\prime}}\subset C_{x,\lambda}$ with $\lambda^{\prime}>\lambda$, then $%
C_{x,\lambda}$ is a neighbourhood of $C_{x^{\prime},\lambda^{\prime}}$.
\end{lemma}

\begin{proof}
$I$ is locally connected, and $h$ is continuous, hence there is a connected
neighbourhood $U$ of $z$ such that $h\left( z\right) \geq\lambda$. By
maximality $U\subset C_{z,\lambda}$. Since $C_{z,\lambda}\cap
C_{y,\lambda}\neq\phi$ we have $C_{z,\lambda}=C_{y,\lambda} $ and thus $%
U\subset C_{y,\lambda}$. Hence $C_{y,\lambda}$ is a neighbourhood of $z$.
The last part follows trivially once by noting that for all $z\in
C_{x^{\prime },\lambda^{\prime}}$ we have $h\left( z\right)
\geq\lambda^{\prime}>\lambda$ and hence $C_{y,\lambda}$ is a neighbourhood
of $z$.
\end{proof}

We now define a pseudo-metric on $I$. Lemma~\ref{lem:loccon} (the only place
we will use local connectedness) is critical to showing that the map from $I$
to the resulting quotient space is continuous.

\begin{definition}
If $y$ and $z$ are points in $I$, define $\lambda\left( y,z\right) \leq
\min\left( h\left( y\right) ,h\left( z\right) \right) $ such that $%
C_{y,\lambda}=C_{z,\lambda}$%
\begin{equation*}
\lambda\left( y,z\right) =\sup\left\{ \lambda\;|\;C_{y,\lambda
}=C_{z,\lambda},\;\lambda\leq h\left( y\right) ,\;\lambda\leq h\left(
z\right) \right\} .
\end{equation*}
\end{definition}

The set 
\begin{equation*}
\left\{ \lambda\;|\;C_{y,\lambda}=C_{z,\lambda},\;\lambda\leq h\left(
y\right) ,\;\lambda\leq h\left( z\right) \right\}
\end{equation*}
is a non-empty interval $\left[ h\left( v\right) ,\lambda\left( y,z\right) %
\right] $ or $\left[ h\left( v\right) ,\lambda\left( y,z\right) \right) $
where $\lambda\left( y,z\right) $ satisfies 
\begin{equation*}
h\left( v\right) \leq\lambda\left( y,z\right) \leq\min\left( h\left(
y\right) ,h\left( z\right) \right) .
\end{equation*}
Clearly $\lambda\left( x,x\right) =h\left( x\right) .$

\begin{lemma}
The function $\lambda $ is lower semi-continuous 
\begin{equation*}
\liminf_{z\rightarrow z_{0}}\lambda \left( y,z\right) \geq \lambda \left(
y,z_{0}\right) .
\end{equation*}
\end{lemma}

\begin{proof}
Fix $y,\;z_{0}$ and choose some $\lambda^{\prime}<\lambda\left(
y,z_{0}\right) $. By the definition of $\lambda(y,z_0)$ we have that $%
C_{y,\lambda^{\prime}}=C_{z_{0},\lambda^{\prime}}$. Since $h\left(
z_{0}\right) \geq\lambda^{\prime}$ there is a neighbourhood $U$ of $z_{0}$
so that $U\subset C_{z_{0},\lambda^{\prime}}$. For any $z\in U$ one has $%
z\in C_{z,\lambda^{\prime}}\cap C_{z_{0},\lambda^{\prime}}$. Hence $%
C_{z_{0},\lambda^{\prime}}=C_{z,\lambda^{\prime}}$ and $C_{y,\lambda^{%
\prime}}=C_{z,\lambda^{\prime}}$. Thus $\lambda\left( y,z\right)
\geq\lambda^{\prime}$ for $z\in U$ and hence 
\begin{equation*}
\liminf_{z\rightarrow z_{0}}\lambda\left( y,z\right) \geq\lambda^{\prime}.
\end{equation*}
Since $\lambda^{\prime}<\lambda\left( y,z_{0}\right) $ was arbitrary 
\begin{equation*}
\liminf_{z\rightarrow z_{0}}\lambda\left( y,z\right) \geq\lambda(y,z_0)
\end{equation*}
and the result is proved.
\end{proof}

\begin{lemma}
The following inequality holds 
\begin{equation*}
\min\left\{ \lambda\left( x,z\right) ,\lambda\left( y,z\right) \right\}
\leq\lambda\left( x,y\right) .
\end{equation*}
\end{lemma}

\begin{proof}
If $\min\left\{ \lambda\left( x,z\right) ,\lambda\left( y,z\right) \right\}
=h\left( v\right)$, then there is nothing to prove. Recall that 
\begin{equation*}
\left\{ \lambda\;|\;C_{y,\lambda}=C_{z,\lambda},\;\lambda\leq h\left(
y\right) ,\;\lambda\leq h\left( x\right) \right\}
\end{equation*}
is connected and contains $h\left( v\right)$. Suppose $h\left( v\right)
\leq\lambda<\min\left\{ \lambda\left( x,z\right) ,\lambda\left( y,z\right)
\right\}$, then it follows that the identity $C_{x,\lambda }=C_{z,\lambda}$
holds for $\lambda$. Similarly $C_{y,\lambda}=C_{z,\lambda}$. As a result $%
C_{x,\lambda}=C_{y,\lambda}$ and $\lambda\left( x,y\right) \geq\lambda$.
\end{proof}

\begin{definition}
Define $d$ on $I\times I$ by 
\begin{equation*}
d\left( x,y\right) =h\left( x\right) +h\left( y\right) -2\lambda\left(
x,y\right).
\end{equation*}
\end{definition}

\begin{lemma}
The function $d$ is a pseudo-metric on $I$. If $\left( \tilde{I},d\right) $
is the resulting quotient metric space, then the projection $I\rightarrow 
\tilde{I}$ from the topological space $I$ to the metric space is continuous.
\end{lemma}

\begin{proof}
Clearly $d$ is positive, symmetric and we have remarked that for all $x$, $%
\lambda\left( x,x\right) =h\left( x\right) $ hence it is zero on the
diagonal. To see the triangle inequality, assume 
\begin{equation*}
\lambda\left( x,z\right) =\min\left\{ \lambda\left( x,z\right)
,\lambda\left( y,z\right) \right\}
\end{equation*}
and then observe 
\begin{align*}
d\left( x,y\right) & =h\left( x\right) +h\left( y\right) -2\lambda\left(
x,y\right) \\
& \leq h\left( x\right) +h\left( y\right) -2\lambda\left( x,z\right) \\
& =h\left( x\right) +h\left( z\right) -2\lambda\left( x,z\right) +h\left(
y\right) -h\left( z\right) \\
& \leq d\left( x,z\right) +\left| h\left( y\right) -h\left( z\right) \right|
\end{align*}
but $\lambda\left( y,z\right) \leq\min\left( h\left( y\right) ,h\left(
z\right) \right) $ and hence 
\begin{align*}
\left| h\left( y\right) -h\left( z\right) \right| & =h\left( y\right)
+h\left( z\right) -2\min\left( h\left( y\right) ,h\left( z\right) \right) \\
& \leq h\left( y\right) +h\left( z\right) -2\lambda\left( y,z\right) \\
& =d\left( y,z\right)
\end{align*}
hence 
\begin{equation*}
d\left( x,y\right) \leq d\left( x,z\right) +d\left( y,z\right)
\end{equation*}
as required.
\end{proof}

We can now introduce the equivalence relation $x\symbol{126}y$ if $d\left(
x,y\right) =0$ and the quotient space $I/\symbol{126}$. We write $I/\symbol{%
126}=\tilde{I}$ and $i:I\rightarrow\tilde{I}$ for the canonical projection.
The function $d$ projects onto $\tilde{I}\times\tilde{I}$ and is a metric
there.

It is tempting to think that $x\symbol{126}y$ if and only if $C_{x}=C_{y}$
and this is true if $I$ is compact Hausdorff. However the definitions imply
a slightly different criteria: $x\symbol{126}y$ iff 
\begin{equation*}
h\left( x\right) =h\left( y\right) =\lambda\text{ and }C_{x,\lambda
^{\prime\prime}}=C_{y,\lambda^{\prime\prime}}\text{ for all }\lambda
^{\prime\prime}<\lambda.
\end{equation*}
The stronger statement $x\symbol{126}y$ if and only if $C_{x}=C_{y}$ is not
true for all continuous functions $h$ on $\mathbb{R}^{2}$ as it is easy to
find a decreasing family of closed connected sets there whose limit is a
closed set that is not connected.

Consider again the new metric space $\tilde{I}$ that has as its points the
equivalence classes of points indistinguishable under $d$. We now prove that
the projection $i$ taking $I$ to $\tilde{I}$ is continuous. Fix $y\in I$ and 
$\varepsilon >0.$ Since $\lambda \left( y,.\right) $ is lower
semi-continuous and $h$ is (upper semi)continuous there is a neighbourhood $%
U $ of $y$ so that for $z\in U$ one has $\lambda \left( y,z\right) >\lambda
\left( y,y\right) -\varepsilon /4$ and $h\left( z\right) <h\left( y\right)
+\varepsilon /2$. Thus $d\left( y,z\right) <\varepsilon $ for $z\in U$.
Hence $\tilde{d}\left( i\left( y\right) ,i\left( z\right) \right)
<\varepsilon $ if $z\in U$. The function $i$ is continuous and as continuous
images of compact sets are compact we have the following.

\begin{corollary}
If $I$ is compact, then $\tilde{I}$ is a compact metric space.
\end{corollary}

To complete this section we will show $\tilde{I}$ is a uniquely arcwise
connected metric space, in which the arc between two points is isometric to
an interval and give a characterisation of compact trees.

\begin{proposition}
If $I$ is a connected and locally connected topological space, and $%
h:I\rightarrow \mathbb{R}$ is a positive continuous function that attains
its lower bound, then its \textquotedblleft contour tree\textquotedblright\
the metric space $\left( \tilde{I},\tilde{d}\right) $ is an $\mathbf{R}$%
-tree. Every $\mathbf{R}$ -tree can be constructed in this way.
\end{proposition}

\begin{proof}
It is enough to prove that the metric space $\tilde{I}$ we have constructed
is really an $\mathbf{R}$-tree and that every $\mathbf{R}$-tree can be
constructed in this way. Let $\tilde{x}$ any point in $\tilde{I}$ and $x\in
I $ satisfy $i\left( x\right) =\tilde{x}$. Then $h\left( x\right) $ does not
depend on the choice of $x$. Fix $h\left( v\right) <\lambda <h\left(
x\right) $. We have seen that there is a $y$ such that $h\left( y\right) $=$%
\lambda $ and $y\prec x$ moreover any two choices have the same contour
through them and hence the same $\tilde{y}\left( \lambda \right) $. In this
way we see that there is a map from $\left[ h\left( v\right) ,h\left(
x\right) \right] $ into $\tilde{I}$ that is injective. Moreover, it is
immediate from the definition of $d$ that it is an isometry and that $\tilde{%
I}$ is uniquely arc connected.

Suppose that $\Omega $ is an $\mathbf{R}$-tree, then we may fix a base
point, and for each point in the tree consider the distance from $V$ it is
clear that this continuous function is just appropriate to ensure that the
contour tree is the original tree.
\end{proof}

\begin{remark}
\textrm{1. In the case where $I$ is compact, obviously $\tilde{I}$ is both
complete and totally bounded as it is compact. }

\textrm{2. An $R$-tree is a metric space; it is therefore possible to
complete it. Indeed the completion consists of those paths, all of whose
initial segments are in the tree\footnote{\textrm{\textrm{We fix a root and
identify the tree with the geodesic arc from the root to the point in the
tree.}}}; we have not identified a simple sufficient condition on the
continuous function and topological space $\Omega $ to ensure this. An $R$%
-tree is totally bounded if it is bounded and for each $\varepsilon >0$
there is an $N$ so that for each $t$ the paths that extend a distance $t$
from the root have at most $N$ ancestral paths between them at time $%
t-\varepsilon $. In this way we see that the $R$-tree that comes out of
studying the historical process for the Fleming-Viot or the Dawson Watanabe
measure-valued processes is, with probability one, a compact $R$-tree for
each finite time.}
\end{remark}

\begin{lemma}
\bigskip Given a compact $R$-tree, there is always a height function on a
closed interval that yields the same tree as its quotient.
\end{lemma}

\begin{proof}
As the tree is compact, path connected and locally path connected, there is
always as based loop mapping $\left[ 0,1\right] $ onto the tree. Let $h$
denote the distance from the root. Its pullback onto the interval $\left[ 0,1%
\right] $ is a height function and the natural quotient is the original
tree. In this way we see that there is always a version of Le Gall's snake 
\cite{snake} traversing a compact tree.
\end{proof}

%\begin{conjecture}
%Given a compact tree $R$-tree it seems reasonable to suppose that there is
%always a continuous function on the closed real line that yields the space
%in the same way as Le Gall's snake traverses the measure valued process.
%\end{conjecture}

%\begin{proof}
%Let $G$ be a group and $\left( E,V\right) $ a connected graph; associate to
%each ordered edge a group element, and to its reverse associate the inverse
%element, then to any path in the graph associate the product of the elements
%associated to the directed edges in the path and in the order they occur. $%
%\left( E,V\right) $ is a tree if and only if for every $G$ and map of $E$
%into $G$, every closed loop in $\left( E,V\right) $ gets mapped to the
%identity element in $G$.
%\end{proof}


\begin{thebibliography}{9}
\bibitem{HL} Hambly, B.M. and Lyons, T.J. \emph{Uniqueness for the signature
of a path of bounded variation and the reduced path group}, to appear Ann.
Math.

\bibitem{Kapovich95} Kapovich, I. \emph{A Non-quasiconvex Subgroup of a
Hyperbolic Group with an Exotic Limit Set}, New York J. Math. 1 (1995),
184-195.

\bibitem{snake} Le Gall, J. F. \emph{Brownian excursions, trees and
measure-valued branching processes} Ann. Probab. 19 (1991), 1399--1439.

\bibitem{rtree} Morgan, J. W., Shalen, P. B., \emph{Valuations, trees, and
degenerations of hyperbolic structures. I.} Ann. of Math. (2) 120 (1984),
401--476.
\end{thebibliography}
\end{document}